\input amstex
\documentstyle{amsppt}
%
\catcode`@=11
\redefine\output@{%
  \def\break{\penalty-\@M}\let\par\endgraf
  \ifodd\pageno\global\hoffset=105pt\else\global\hoffset=8pt\fi  
  \shipout\vbox{%
    \ifplain@
      \let\makeheadline\relax \let\makefootline\relax
    \else
      \iffirstpage@ \global\firstpage@false
        \let\rightheadline\frheadline
        \let\leftheadline\flheadline
      \else
        \ifrunheads@ 
        \else \let\makeheadline\relax
        \fi
      \fi
    \fi
    \makeheadline \pagebody \makefootline}%
  \advancepageno \ifnum\outputpenalty>-\@MM\else\dosupereject\fi
}
\def\Beta{\mathchar"0\hexnumber@\rmfam 42}
\catcode`\@=\active
\nopagenumbers
\chardef\textvolna='176

\chardef\bigalpha='013
\def\negskp{\hskip -2pt}
\def\Rea{\operatorname{Re}}
\def\Img{\operatorname{Im}}

\chardef\degree="5E
\def\blue#1{#1}

\gdef\darkred#1{#1}
\catcode`#=11\def\diez{#}\catcode`#=6
\catcode`&=11\catcode`&=4
\catcode`_=11\def\podcherkivanie{_}\catcode`_=8
\catcode`\^=11\def\shlyapka{^}\catcode`\^=7
\catcode`~=11\def\volna{~}\catcode`~=\active
\def\mycite#1{\cite{\blue{#1}}\immediate\special{ps:
     ShrHPSdict begin /ShrBORDERthickness 0 def}}
\def\myciterange#1#2#3#4{\cite{\blue{#2#3#4}}\immediate\special{ps:
     ShrHPSdict begin /ShrBORDERthickness 0 def}}
\def\mytag#1{%
    \tag#1}
\def\mythetag#1{\thetag{\blue{#1}}\immediate\special{ps:
     ShrHPSdict begin /ShrBORDERthickness 0 def}}
\def\myrefno#1{\no#1}
\def\myhref#1#2{\blue{#2}\immediate\special{ps:
     ShrHPSdict begin /ShrBORDERthickness 0 def}}
\def\myEarXivlink{\myhref{http://arXiv.org}{http:/\negskp/arXiv.org}}

\def\mytheorem#1{\csname proclaim\endcsname{Theorem #1}}
\def\mytheoremwithtitle#1#2{\csname proclaim\endcsname{Theorem #1#2}}
\def\mythetheorem#1{\blue{#1}\immediate\special{ps:
     ShrHPSdict begin /ShrBORDERthickness 0 def}}
\def\mylemma#1{\csname proclaim\endcsname{Lemma #1}}
\def\mylemmawithtitle#1#2{\csname proclaim\endcsname{Lemma #1#2}}
\def\mythelemma#1{\blue{#1}\immediate\special{ps:
     ShrHPSdict begin /ShrBORDERthickness 0 def}}
\def\mycorollary#1{\csname proclaim\endcsname{Corollary #1}}

\def\mydefinition#1{\definition{Definition #1}}
\def\mythedefinition#1{\blue{#1}\immediate\special{ps:
     ShrHPSdict begin /ShrBORDERthickness 0 def}}
\def\myconjecture#1{\csname proclaim\endcsname{Conjecture #1}}
\def\myconjecturewithtitle#1#2{\csname proclaim\endcsname{Conjecture #1#2}}
\def\mytheconjecture#1{\blue{#1}\immediate\special{ps:
     ShrHPSdict begin /ShrBORDERthickness 0 def}}
\def\myproblem#1{\csname proclaim\endcsname{Problem #1}}
\def\myproblemwithtitle#1#2{\csname proclaim\endcsname{Problem #1#2}}


\pagewidth{360pt}
\pageheight{606pt}
\topmatter
\title
A strategy of numeric search for perfect cuboids in the case
of the second cuboid conjecture.
\endtitle
\rightheadtext{A strategy of numeric search for perfect cuboids \dots}
\author
A.\,A.\,Masharov, R.\,A.\,Sharipov
\endauthor
\address Bashkir State University, 32 Zaki Validi street, 450074 Ufa, Russia
\endaddress
\email\vtop{\hsize=6cm\noindent\myhref{mailto:masharov92\@mail.ru}
{masharov92\@mail.ru}\newline
\myhref{mailto:r-sharipov\@mail.ru}{r-sharipov\@mail.ru}}
\endemail
\abstract
     A perfect cuboid is a rectangular parallelepiped whose edges, whose face 
diagonals, and whose space diagonal are of integer lengths. The problem of 
finding such cuboids or proving their non-existence is not solved thus far. 
The second cuboid conjecture specifies a subclass of perfect cuboids described
by one Diophantine equation of tenth degree and claims their non-existence
within this subclass. Regardless of proving or disproving this conjecture 
in the present paper the Diophantine equation associated with it is studied
and is used in order to build an optimized strategy of computer-assisted search 
for perfect cuboids within the subclass covered by the second cuboid conjecture.   
\endabstract
\endtopmatter
\TagsOnRight
\document

\head
1. Introduction.
\endhead
     For the history and various approaches to the problem of perfect 
cuboids the reader is referred to \myciterange{1}{1}{--}{42}. In this 
paper we resume the research initiated in \myciterange{43}{43}{--}{47}. 
The papers \myciterange{48}{48}{--}{60} deal with another approach based on 
so-called multisymmetric polynomials. In this paper we do not touch this 
approach.\par
     Perfect cuboids are described by six Diophantine equations. These
equations are immediate from the Pythagorean theorem:  
$$
\xalignat 2
&\hskip -2em
x_1^2+x_2^2+x_3^2-L^2=0,
&&x_2^2+x_3^2-d_1^{\kern 1pt 2}=0,\\
\vspace{-1.7ex}
\mytag{1.1}\\
\vspace{-1.7ex}
&\hskip -2em
x_3^2+x_1^2-d_2^{\kern 1pt 2}=0,
&&x_1^2+x_2^2-d_3^{\kern 1pt 2}=0.
\endxalignat
$$
The variables $x_1$, $x_2$, $x_3$ in \mythetag{1.1} stand for three edges 
of a cuboid, the variables $d_1$, $d_2$, $d_3$ correspond to its face diagonals, 
and $L$ represents its space diagonal.\par
      In \mycite{43} an algebraic parametrization for the Diophantine equations 
\mythetag{1.1} was suggested. It uses four rational variables 
$\alpha$, $\beta$, $\upsilon$, and $z$:
$$
\xalignat 2
&\hskip -2em 
\frac{x_1}{L}=\frac{2\,\upsilon}{1+\upsilon^2},
&&\frac{d_1}{L}=\frac{1-\upsilon^2}{1+\upsilon^2},\\
\vspace{1ex}
&\hskip -2em
\frac{x_2}{L}=\frac{2\,z\,(1-\upsilon^2)}{(1+\upsilon^2)\,(1+z^2)},
&&\frac{x_3}{L}=\frac{(1-\upsilon^2)\,(1-z^2)}{(1+\upsilon^2)\,(1+z^2)},
\mytag{1.2}\\
\vspace{1ex}
&\hskip -2em
\frac{d_2}{L}=\frac{(1+\upsilon^2)\,(1+z^2)+2\,z(1-\upsilon^2)}
{(1+\upsilon^2)\,(1+z^2)}\,\beta,
&&\frac{d_3}{L}=\frac{2\,(\upsilon^2\,z^2+1)}{(1+\upsilon^2)\,(1+z^2)}\,\alpha.
\quad
\endxalignat
$$
The variables $\alpha$, $\beta$, $\upsilon$ in \mythetag{1.2} are different 
from the original ones which are used in \mycite{43}, here we use 
$\alpha$ and $\beta$ instead of $a$ and $b$, and we use $\upsilon$ instead 
of $u$.\par 
    Only two of the four variables $\alpha$, $\beta$, $\upsilon$, and $z$ are 
independent. The variables $\alpha$ and $\beta$ are taken for independent 
ones. Then the variable $\upsilon$ is expressed through $\alpha$ and $\beta$ 
as a solution of the following algebraic equation:
$$
\gathered
\upsilon^4\,\alpha^4\,\beta^4+(6\,\alpha^4\,\upsilon^2\,\beta^4
-2\,\upsilon^4\,\alpha^4\,\beta^2-2\,\upsilon^4\,\alpha^2
\,\beta^4)+(4\,\upsilon^2\,\beta^4\,\alpha^2+\\
+\,4\,\alpha^4\,\upsilon^2\,\beta^2-12\,\upsilon^4\,\alpha^2\,\beta^2
+\upsilon^4\,\alpha^4+\upsilon^4\,\beta^4
+\alpha^4\,\beta^4)+(6\,\alpha^4\,\upsilon^2+6\,\upsilon^2\,\beta^4-\\
-\,8\,\alpha^2\,\beta^2\,\upsilon^2-2\,\upsilon^4\,\alpha^2
-2\,\upsilon^4\,\beta^2-2\,\alpha^4\,\beta^2-2\,\beta^4\,\alpha^2)
+(\upsilon^4+\beta^4+\\
+\,\alpha^4+4\,\alpha^2\,\upsilon^2+4\,\beta^2\,\upsilon^2
-12\,\beta^2\,\alpha^2)+(6\,\upsilon^2-2\,\alpha^2-2\,\beta^2)+1=0.
\endgathered
\quad
\mytag{1.3}
$$
Once the variable $\upsilon$ is expressed as a function of $\alpha$ and 
$\beta$ by solving the equation \mythetag{1.3}, the variable $z$ is
given by the formula 
$$
\hskip -2em
z=\frac{(1+\upsilon^2)\,(1-\beta^2)\,(1+\alpha^2)}{2\,(1+\beta^2)\,(1
-\alpha^2\,\upsilon^2)}.
\mytag{1.4}
$$
The equation \mythetag{1.3}, along with the formula \mythetag{1.4}, 
produces two algebraic functions
$$
\xalignat 2
&\hskip -2em
\upsilon=\upsilon(\alpha,\beta),
&&z=z(\alpha,\beta).
\mytag{1.5}
\endxalignat
$$
Substituting \mythetag{1.5} into \mythetag{1.2}, we get six algebraic
functions
$$
\xalignat 2
&\hskip -2em
x_1=x_1(\alpha,\beta,L),
&&d_1=d_1(\alpha,\beta,L),\\
&\hskip -2em
x_2=x_2(\alpha,\beta,L),
&&d_2=d_2(\alpha,\beta,L),
\mytag{1.6}\\
&\hskip -2em
x_3=x_3(\alpha,\beta,L),
&&d_3=d_3(\alpha,\beta,L),
\endxalignat
$$
which are linear with respect to $L$. The functions \mythetag{1.6} satisfy 
the cuboid equations \mythetag{1.1} identically with respect to $\alpha$, 
$\beta$, and $L$. This fact is presented by the following theorem 
(see Theorem~5.2 in \mycite{43}). 
\mytheorem{1.1} A perfect cuboid does exist if and only if there are
three rational numbers $\alpha$, $\beta$, and $\upsilon$ satisfying the 
equation \mythetag{1.3} and obeying four inequalities $0<\alpha<1$,
$0<\beta<1$, $0<\upsilon<1$, and $(\alpha+1)\,(\beta+1)>2$.
\endproclaim
     The rational numbers $\alpha$, $\beta$, and $\upsilon$ cam be brought 
to a common denominator:
$$
\xalignat 3
&\hskip -2em
\alpha=\frac{a}{t},
&&\beta=\frac{b}{t},
&&\upsilon=\frac{u}{t}.
\mytag{1.7}
\endxalignat
$$
Substituting \mythetag{1.7} into \mythetag{1.3}, one easily derives the
Diophantine equation
$$
\gathered
t^{12}+(6\,u^2\,-2\,a^2\,-2\,b^2)\,t^{10}
+(u^4\,+b^4+a^4+4\,a^2\,u^2+\\
+\,4\,b^2\,u^2-12\,b^2\,a^2)\,t^8+(6\,a^4\,u^2+6\,u^2\,b^4-8\,a^2\,b^2\,u^2-\\
-\,2\,u^4\,a^2-2\,u^4\,b^2-2\,a^4\,b^2-2\,b^4\,a^2)\,t^6+(4\,u^2\,b^4\,a^2+\\
+\,4\,a^4\,u^2\,b^2-12\,u^4\,a^2\,b^2+u^4\,a^4+u^4\,b^4+a^4\,b^4)\,t^4+\\
+\,(6\,a^4\,u^2\,b^4-2\,u^4\,a^4\,b^2-2\,u^4\,a^2
\,b^4)\,t^2+u^4\,a^4\,b^4=0.
\endgathered
\quad
\mytag{1.8}
$$
Theorem~\mythetheorem{1.1} then is reformulated in the following form
(see Theorem 4.1 in \mycite{44}). 
\mytheorem{1.2} A perfect cuboid does exist if and only if for some
positive coprime integer numbers $a$, $b$, and $u$ the Diophantine
equation \mythetag{1.8} has a positive solution $t$ obeying the inequalities
$t>a$, $t>b$, $t>u$, and $(a+t)\,(b+t)>2\,t^2$. 
\endproclaim
     In \mycite{44} the Diophantine equation was treaded as a polynomial 
equation for $t$, while $a$, $b$, and $u$ were considered as parameters. 
As a result in \mycite{44} several special cases of the equation 
\mythetag{1.8} were specified. They are introduced through the following 
relationships for the parameters $a$, $b$, and $u$:
$$
\xalignat 3
&\hskip -2em
\text{1) \ }a=b\neq u; &&\text{3) \ }b\,u=a^2; &&\text{5) \ }a=u\neq b;\\
\vspace{-1.5ex}
\mytag{1.9}\\
\vspace{-1.5ex}
&\hskip -2em
\text{2) \ }a=b=u; &&\text{4) \ }a\,u=b^2; &&\text{6) \ }b=u\neq a.
\endxalignat
$$
The cases 2, 5, and 6 are trivial. They produce no perfect cuboids
(see \mycite{44}). The case 1 corresponds to the first cuboid conjecture
(see \mycite{44}). It is less trivial, but it produces no perfect 
cuboids either (see \mycite{45}). The cases 2 and 4 correspond to the 
second cuboid conjecture (see \mycite{44} and \mycite{46}). The case, 
where none of the conditions \mythetag{1.9} is fulfilled, corresponds 
to the third cuboid conjecture (see \mycite{44} and \mycite{47}).\par
     In this paper we consider the cases 3 and 4 associated with the 
second cuboid conjecture. In the case 3 the equality $b\,u=a^2$ is
resolved by substituting 
$$
\xalignat 3
&\hskip -2em
a=p\,q,
&&b=p^{\kern 1pt 2}, &&u=q^{\kern 0.7pt 2}.
\mytag{1.10}
\endxalignat
$$
Here $p\neq q$ are two positive coprime integers. Upon substituting 
\mythetag{1.10} into the equation \mythetag{1.8} it reduces to the 
equation 
$$
(t-a)\,(t+a)\,Q_{pq}(t)=0
\mytag{1.11}
$$
(see \mycite{45}), where $Q_{pq}(t)$ is the following polynomial of tenth
degree:
$$
\gathered
Q_{pq}(t)=t^{10}+(2\,q^{\kern 0.7pt 2}+p^{\kern 1pt 2})\,(3\,q^{\kern 0.7pt 2}-2\,p^{\kern 1pt 2})\,t^8
+(q^{\kern 0.5pt 8}+10\,p^{\kern 1pt 2}\,q^{\kern 0.5pt 6}+\\
+\,4\,p^{\kern 1pt 4}\,q^4-14\,p^{\kern 1pt 6}\,q^{\kern 0.7pt 2}+p^{\kern 1pt 8})\,t^6
-p^{\kern 1pt 2}\,q^{\kern 0.7pt 2}\,(q^{\kern 0.5pt 8}-14\,p^{\kern 1pt 2}\,q^{\kern 0.5pt 6}+4\,p^{\kern 1pt 4}\,q^4+\\
+\,10\,p^{\kern 1pt 6}\,q^{\kern 0.7pt 2}+p^{\kern 1pt 8})\,t^4
-p^{\kern 1pt 6}\,q^{\kern 0.5pt 6}\,(q^{\kern 0.7pt 2}
+2\,p^{\kern 1pt 2})\,(3\,p^{\kern 1pt 2}-2\,q^{\kern 0.7pt 2})\,t^2
-q^{\kern 0.7pt 10}\,p^{\kern 1pt 10}.
\endgathered
\quad
\mytag{1.12}
$$\par
     The case 4 is similar. In this case the equality $a\,u=b^2$ is
resolved by substituting 
$$
\xalignat 3
&\hskip -2em
a=p^{\kern 1pt 2},
&&b=p\,q, &&u=q^{\kern 0.7pt 2}.
\mytag{1.13}
\endxalignat
$$
Upon substituting \mythetag{1.13} into the equation \mythetag{1.8} 
it reduces to the equation 
$$
\hskip -2em
(t-b)\,(t+b)\,Q_{pq}(t)=0.
\mytag{1.14}
$$
The roots $t=a$, $t=-a$, $t=b$, and $t=-b$ of the equations \mythetag{1.11}
and \mythetag{1.14} do not produce perfect cuboids 
(see Theorem~\mythetheorem{1.2}). Upon splitting off the linear factors
from \mythetag{1.11} and \mythetag{1.14} we get the equation 
$$
\gathered
t^{10}+(2\,q^{\kern 0.7pt 2}+p^{\kern 1pt 2})\,(3\,q^{\kern 0.7pt 2}
-2\,p^{\kern 1pt 2})\,t^8+(q^{\kern 0.5pt 8}+10\,p^{\kern 1pt 2}\,q^{\kern 0.5pt 6}
+4\,p^{\kern 1pt 4}\,q^4\,-\\
-\,14\,p^{\kern 1pt 6}\,q^{\kern 0.7pt 2}+p^{\kern 1pt 8})\,t^6
-p^{\kern 1pt 2}\,q^{\kern 0.7pt 2}\,(q^{\kern 0.5pt 8}
-14\,p^{\kern 1pt 2}\,q^{\kern 0.5pt 6}+4\,p^{\kern 1pt 4}\,q^4
+10\,p^{\kern 1pt 6}\,q^{\kern 0.7pt 2}+\\
+p^{\kern 1pt 8})\,t^4-p^{\kern 1pt 6}\,q^{\kern 0.5pt 6}\,(q^{\kern 0.7pt 2}
+2\,p^{\kern 1pt 2})\,(3\,p^{\kern 1pt 2}-2\,q^{\kern 0.7pt 2})\,t^2
-q^{\kern 0.7pt 10}\,p^{\kern 1pt 10}=0.
\endgathered
\quad
\mytag{1.15}
$$
\myconjecture{1.1} For any positive coprime integers $p\neq q$ the 
polynomial $Q_{pq}(t)$ in \mythetag{1.12} is irreducible in the ring 
$\Bbb Z[t]$. 
\endproclaim
     Conjecture~\mytheconjecture{1.1} is known as the second cuboid
conjecture. It was formulated in \mycite{44}. In particular it claims 
that the equation \mythetag{1.15} has no integer roots for any positive 
coprime integers $p\neq q$. We do not try to prove or disprove 
Conjecture~\mytheconjecture{1.1} in this paper. Instead, we study
real positive roots of the equation \mythetag{1.15} in the case where
$q$ is much larger than $p$. Using asymptotic expansions for the roots
of the equation \mythetag{1.15} as $q\to+\infty$, below we build an 
optimized strategy of computer-assisted search for perfect cuboids in 
the realm of Conjecture~\mytheconjecture{1.1}.\par
\head
2. Asymptotic expansions for roots of the polynomial equation.
\endhead
     Note that the polynomial $Q_{pq}(t)$ in \mythetag{1.12} is even. 
Along with each root $t$ it has the opposite root $-t$. We use the 
condition
$$
\hskip -2em
\cases \text{$t>0$ \ if \ $t$ \ is a real root,}\\
\text{$\Rea(t)\geqslant 0$ \ and \ $\Img(t)>0$ \ if \ $t$ \ is a complex root,}
\endcases
\mytag{2.1}
$$
in order to divide the roots of the equation \mythetag{1.15} into
two groups. We denote through $t_1$, $t_2$, $t_3$, $t_4$, $t_5$ the
roots that obey the conditions \mythetag{2.1}. Then $t_6$, $t_7$, $t_8$, 
$t_9$, $t_{10}$ are opposite roots of the equation \mythetag{1.15}:
$$
\xalignat 5
&\hskip -1em
t_6=-t_1,&&t_7=-t_2,&&t_8=-t_3,&&t_9=-t_4,&&t_{10}=-t_5.
\qquad\quad
\mytag{2.2}
\endxalignat
$$\par
     Typically, asymptotic expansions for roots of a polynomial equation 
look like power series (see \mycite{61}). In our case we have the expansions
$$
\hskip -2em
t_i(p,q)=C_i\,q^{\,\alpha_i}\biggl(1+\sum^\infty_{s=1}\beta_{is}\,q^{-s}\biggr)
\text{\ \ as \ }q\to+\infty.
\mytag{2.3}
$$
The coefficient $C_i$ in \mythetag{2.3} should be nonzero: $C_i\neq 0$.\par
     Let's substitute \mythetag{2.3} into the equation \mythetag{1.15}. For this
purpose we represent the polynomial $Q_{pq}(t)$ from \mythetag{1.12} 
formally as the sum 
$$
\hskip -2em
Q_{pq}(t)=\sum^{10}_{m=0}\sum^{10}_{r=0}A_{m\kern 1pt r}(p)\,q^{\kern 1pt r}
\,t^m.
\mytag{2.4}
$$
Each nonzero term in \mythetag{2.4}, i\.\,e\. a term $A_{m\kern 1pt r}(p)
\,q^{\kern 1pt r}\,t^m$ with the nonzero coefficient 
$$
\gather
\hskip -2em
A_{m\kern 1pt r}(p)\neq 0,\\
\vspace{-2ex}
\intertext{yields the sum}
\hskip -2em
S_{m\kern 1pt r}(p,q)=A_{m\kern 1pt r}(p)\ {C_i}^m
\,q^{\kern 1pt m\,\alpha_i+r}\kern 0.6em
+\kern -0.6em\sum_{s<m\,\alpha_i+r}\kern -1em\gamma_{irms}\,q^s.
\mytag{2.5}
\endgather
$$
Taking into account \mythetag{2.5}, the equation \mythetag{1.15} 
is written as
$$
\hskip -2em
\sum^{10}_{m=0}\sum^{10}_{r=0}S_{m\kern 1pt r}(p,q)=0.
\mytag{2.6}
$$
\par
     The equality \mythetag{2.6} should be fulfilled identically with respect 
to the variable $q\to+\infty$. Since $C_i\neq 0$, a necessary condition for 
that is the coincidence of exponents of at least two summands of the form
$S_{m\kern 1pt r}(q)$ in the leading order with respect to the variable $q$.
This yields the equalities 
$$
\hskip -2em
m_1\,\alpha_i+r_1=m_2\,\alpha_i+r_2=s_{\sssize\text{max}}. 
\mytag{2.7}
$$ 
The maximality of the exponent in \mythetag{2.7} means that all exponents
are not greater than $s_{\sssize\text{max}}$, i\.\,e\. we have the following
inequality: 
$$
\hskip -2em
m\,\alpha_i+r\leqslant s_{\sssize\text{max}}
\text{\ \ for all \ $r$ \ and \ $m$\ \ such that \ $A_{m\kern 1pt r}(p)\neq 0$.}
\mytag{2.8}
$$
\mylemma{2.1} The coincidence $m_1=m_2$ in the formula \mythetag{2.7} 
is impossible. 
\endproclaim 
\demo{Proof} Indeed, due to \mythetag{2.8} the coincidence $m_1=m_2$ would mean
$r_1=r_2$. But the sum \mythetag{2.6} has no two summands with simultaneously 
coinciding indices $r$ and $m$. Lemma~\mythelemma{2.1} is proved.\qed\enddemo
     Let's treat $m$ and $r$ as coordinates of a point on the coordinate plane. 
Since $m$ and $r$ are integer, such a point belongs to the integer grid, being
its node. \vadjust{\vskip 270pt\hbox to 0pt{\kern 15pt
\includegraphics{Strategy01.eps}\hss}\vskip 0pt}The numbers $m_1$, 
$r_1$ and $m_2$, $r_2$ from \mythetag{2.7} mark two nodes  of this grid. 
These are the points $A$ and $B$ in Fig\.~2.1. Due to Lemma~\mythelemma{2.1} 
from the equality \mythetag{2.7} we derive the following formula for the 
exponent $\alpha_i$:
$$
\hskip -2em
\alpha_i=-\frac{r_2-r_1}{m_2-m_1}.
\mytag{2.9}
$$ 
The right hand side of the formula \mythetag{2.9} up to the sign coincides 
with the slope of the straight line connecting the nodes $A$ and $B$ in
Fig\.~2.1:
$$
\hskip -2em
\alpha_i=-k_{\sssize AB}.
\mytag{2.10}
$$\par
     The nodes $A$ and $B$ correspond to some nonzero summands in the sum
\mythetag{2.4} being a formal presentation of the polynomial \mythetag{1.12}. 
They are selected by the maximality condition for the parameter $s=m\,\alpha_i+r$. 
The maximum is taken over all summands in the sum \mythetag{2.4} for a fixed value
of $\alpha_i$. 
\mylemma{2.2} The exponent $\alpha_i$ in the asymptotic expansion \mythetag{2.3} 
is determined by the slope of a straight line connecting some two nodes of the 
integer grid associated with some two nonzero terms in the polynomial 
\mythetag{2.4}. 
\endproclaim
    Let $C$ be some node of the integer grid in Fig\.~2.1 associated with some
nonzero summand of the sum \mythetag{2.4} and different from the nodes $A$ and
$B$. Its coordinates $m$ and $r$ satisfy the inequality \mythetag{2.8}. From
\mythetag{2.7} and \mythetag{2.8} one derives the inequality 
$$
m\,\alpha_i+r\leqslant m_1\,\alpha_i+r_1.
$$
Let's write this inequality as follows:
$$
\hskip -2em
\alpha_i\,(m-m_1)\leqslant -(r-r_1).
\mytag{2.11}
$$\par
     In Fig\.~3.1 three positions of the node $C$ relative to the node $A$ are
shown. The node $C$ can be located to the left of the node $A$, to the right of
the node $A$, or on the same vertical line with the node $A$. In the first case 
$m<m_1$. In the second case $m>m_1$. And finally, in the third case $m=m_1$.\par
     In the first case, i\.\,e\. if $m<m_1$, from \mythetag{2.11} we derive 
$$
\hskip -2em
\alpha_i\geqslant -\frac{r-r_1}{m-m_1}.
\mytag{2.12}
$$
The right hand side of the inequality \mythetag{2.12} up to the sign coincides 
with the slope of the line $AC$. Applying \mythetag{2.10}, we get the inequality 
$-k_{\ssize AB}\geqslant -k_{\ssize AC}$. Inverting signs, we write this inequality 
in the following form: 
$$
\hskip -2em
k_{\ssize AC}\geqslant k_{\ssize AB}.
\mytag{2.13}
$$\par
     In the second case, i\.\,e\. if $m>m_1$, from \mythetag{2.11} we derive 
$$
\hskip -2em
\alpha_i\leqslant -\frac{r-r_1}{m-m_1}.
\mytag{2.14}
$$
By analogy with \mythetag{2.13} the inequality \mythetag{2.14} is transformed 
to 
$$
\hskip -2em
k_{\ssize AC}\leqslant k_{\ssize AB}.
\mytag{2.15}
$$\par
     And finally, in the third case, i\.\,e\. if  $m=m_1$, from the inequality 
\mythetag{2.11} we derive 
$$
\hskip -2em
0\leqslant -(r-r_1).
\mytag{2.16}
$$
The inequality \mythetag{2.16} is equivalent to the following inequality:
$$
\pagebreak 
\hskip -2em
r\leqslant r_1.
\mytag{2.17}
$$
Each of the inequalities \mythetag{2.13}, \mythetag{2.15}, and  \mythetag{2.17} 
in its case means that the point $C$ is located not above the line $AB$. This fact
is formulated as a lemma.
\mylemma{2.3} All nodes $(m,r)$ of the integer grid associated with nonzero 
summands in the polynomial \mythetag{2.4} are located not above the line $AB$
on which the nodes implementing the maximum of the parameter $s=m\,\alpha_i+r$ 
are located. 
\endproclaim
     In order to apply Lemmas~\mythelemma{2.1}, \mythelemma{2.2}, and 
\mythelemma{2.3} let's mark all of the nodes associated with the polynomial
\mythetag{1.12} \vadjust{\vskip 364pt\hbox to 0pt{\kern 10pt 
\includegraphics{Strategy02.eps}\hss}\vskip -4pt}on the coordinate 
plane.\par
\mydefinition{2.1} For any polynomial of two variables $P(t,q)$ the convex 
hull of all integer nodes on the coordinate plane associated with monomials 
of this polynomial is called the Newton polygon of $P(t,q)$.  
\enddefinition
{\bf Remark}. Note that in our case the polynomial \mythetag{1.12} depend 
on three variables $p$, $q$, and $t$. However, we treat $p$ as a parameter 
and consider $Q_{pq}(t)$ as a polynomial of two variables when applying 
Definition~\mythedefinition{2.1} to it.\par
      The Newton polygon of the polynomial \mythetag{1.12} is shown in 
Fig\.~2.2. Its boundary consists of two parts --- the upper part and the 
lower part. The upper parts is drawn in green, the lower part is drawn in red. 
In Fig\.~2.2 the nodes on the upper boundary of the Newton polygon are denoted 
according to the formula \mythetag{2.4}. The coefficients $A_{m\kern 1pt r}(p)$
in \mythetag{1.12} associated with these nodes are given by the formulas
$$
\xalignat 3
&\hskip -2em
A_{\kern 0.5pt 0\kern 2pt 10}=-p^{\kern 1pt 10},
&&A_{\kern 0.5pt 2\kern 1.5pt 10}=2\,p^{\kern 1pt 6},
&&A_{\kern 0.5pt 4\kern 1.5pt 10}=-p^{\kern 1pt 2},
\quad
\\
\vspace{-1.5ex}
\mytag{2.18}\\
\vspace{-1.5ex}
&\hskip -2em
A_{\kern 0.5pt 6\kern 2pt 8}=1,
&&A_{\kern 0.5pt 8\kern 1.5pt 4}=6,
&&A_{\kern 0.5pt 10\kern 2pt 0}=1.
\quad
\endxalignat
$$
\mytheorem{2.1} The values of exponents $\alpha_i$ in the expansion \mythetag{2.3} 
for roots of the equation \mythetag{1.15} are determined according to the formula
$\alpha_i=-k$, where $k$ stands for slopes of segments of the polygonal line being 
the upper boundary of the Newton polygon in Fig\.~2.2. 
\endproclaim
     Theorem~\mythetheorem{2.1} is immediate from Lemmas~\mythelemma{2.2} and
\mythelemma{2.3}. The formula $\alpha_i=-k$ in this theorem follows from the formula~\mythetag{2.10}. In our particular case we have
$$
\xalignat 3
&\hskip -2em
\alpha_i=0,
&&\alpha_i=1,
&&\alpha_i=2.
\mytag{2.19}
\endxalignat
$$
The options \mythetag{2.19} are derived from Fig\.~2.2 due to the above theorem.
\par
\head
3. Leading terms in asymptotic expansions.
\endhead
     The term $C_i\,q^{\,\alpha_i}$ obtained upon expanding brackets in 
\mythetag{2.3} is called the leading term of the asymptotic expansion 
\mythetag{2.3}. Three options for the value of $\alpha_i$ are given by the
formula \mythetag{2.19}. Let's consider each of these options separately. 
\par
    {\bf The case }$\alpha_i=0$. This case corresponds to the horizontal 
segment on the upper boundary of the Newton polygon in Fig\.~2.2. This 
segment comprises three nodes $A_{\kern 0.5pt 4\kern 1.5pt 10}$,
$A_{\kern 0.5pt 2\kern 1.5pt 10}$, and $A_{\kern 0.5pt 0\kern 2pt 10}$. 
Therefore, substituting the expansion \mythetag{2.3} with $\alpha_i=0$ 
into the equation \mythetag{1.15}, we get the following equation for
$C_i$:
$$
\hskip -2em
A_{\kern 0.5pt 4\kern 1.5pt 10}\ {C_i}^4
+A_{\kern 0.5pt 2\kern 1.5pt 10}\ {C_i}^2
+A_{\kern 0.5pt 0\kern 2pt 10}=0. 
\mytag{3.1}
$$
Taking into account \mythetag{2.18}, the equation \mythetag{3.1} is 
transformed to
$$
\hskip -2em
p^{\kern 1pt 2}\,{C_i}^4-2\,p^{\kern 1pt 6}\,{C_i}^2+p^{\kern 1pt 10}=0.
\mytag{3.2}
$$
The equation \mythetag{3.2} has two real roots 
$$
\xalignat 2
&\hskip -2em
C_i=p^{\kern 1pt 2},
&&C_i=-p^{\kern 1pt 2},
\mytag{3.3}
\endxalignat
$$
each of which is of multiplicity $2$. The condition \mythetag{2.1} excludes
the root $C_i=-1$ from \mythetag{3.3}. The remain is one root of multiplicity $2$:
$$
\hskip -2em
C_i=p^{\kern 1pt 2}.
\mytag{3.4}
$$
The asymptotic expansion \mythetag{2.3} corresponding to \mythetag{3.4} is
$$
\hskip -2em
t_i(p,q)=p^{\kern 1pt 2}\biggl(1+\sum^\infty_{s=1}\beta_{is}\,q^{-s}\biggr). 
\mytag{3.5}
$$\par
     {\bf The case} $\alpha_i=1$. This case corresponds to the short slant segment 
in the upper boundary of the Newton polygon in Fig\.~2.2. It comprises
two nodes $A_{\kern 0.5pt 4\kern 1.5pt 10}$ and $A_{\kern 0.5pt 6\kern 2pt 8}$. 
Therefore, substituting the expansion \mythetag{2.3} with $\alpha_i=1$ 
into the equation \mythetag{1.15}, we get the following equation for $C_i$:
$$
\hskip -2em
A_{\kern 0.5pt 6\kern 2pt 8}\ {C_i}^6+
A_{\kern 0.5pt 4\kern 1.5pt 10}\ {C_i}^4=0.
\mytag{3.6}
$$
The common divisor ${C_i}^4$ can be factored out from the equation \mythetag{3.6}.
Since $C_i\neq 0$, we can remove this common divisor. Then the equation takes the 
form 
$$
\hskip -2em
A_{\kern 0.5pt 6\kern 2pt 8}\ {C_i}^2+
A_{\kern 0.5pt 4\kern 1.5pt 10}=0.
\mytag{3.7}
$$
Taking into account \mythetag{2.18}, the equation \mythetag{3.7} is
transformed to
$$
\hskip -2em
{C_i}^2-p^{\kern 1pt 2}=0.
\mytag{3.8}
$$
The quadratic equation \mythetag{3.8} has two simple root
$$
\xalignat 2
&\hskip -2em
C_i=p,
&&C_i=-p,
\mytag{3.9}
\endxalignat
$$
The condition \mythetag{2.1} excludes the root $C_i=-p$ from \mythetag{3.9}. 
Therefore as a remain we have only one root, which is of multiplicity $1$:  
$$
\hskip -2em
C_i=p.
\mytag{3.10}
$$
The asymptotic expansion \mythetag{2.3} corresponding to \mythetag{3.10} is
$$
\hskip -2em
t_i(p,q)=p\,q\biggl(1+\sum^\infty_{s=1}\beta_{is}\,q^{-s}\biggr).
\mytag{3.11}
$$\par
     {\bf The case} $\alpha_i=2$. This case corresponds to the long slant 
segment in the upper boundary of the Newton polygon in Fig\.~3.2. It 
comprises three nodes $A_{\kern 0.5pt 6\kern 2pt 8}$, 
$A_{\kern 0.5pt 8\kern 1.5pt 4}$, and $A_{\kern 0.5pt 10\kern 2pt 0}$.
Therefore, substituting the expansion \mythetag{2.3} with $\alpha_i=2$ 
into the equation \mythetag{1.15}, we get the following equation for $C_i$:     
$$
\hskip -2em
A_{\kern 0.5pt 10\kern 2pt 0}\ {C_i}^{10}
+A_{\kern 0.5pt 8\kern 1.5pt 4}\ {C_i}^8
+A_{\kern 0.5pt 6\kern 2pt 8}\ {C_i}^6=0.
\mytag{3.12}
$$
The common divisor ${C_i}^6$ is factored out from the equation \mythetag{3.12}.
Since $C_i\neq 0$, we can remove this common divisor. Then the equation takes the 
form 
$$
\hskip -2em
A_{\kern 0.5pt 10\kern 2pt 0}\ {C_i}^4
+A_{\kern 0.5pt 8\kern 1.5pt 4}\ {C_i}^2
+A_{\kern 0.5pt 6\kern 2pt 8}=0.
\mytag{3.13}
$$
Taking into account \mythetag{2.18}, the equation \mythetag{3.13} is
transformed to
$$
\hskip -2em
{C_i}^4+6\ {C_i}^2+1=0.
\mytag{3.14}
$$
The quartic equation \mythetag{3.14} has four roots. All of them are complex:
$$
\xalignat 2
&\hskip -2em
C_i=(\sqrt{2}+1)\,\goth i,
&&C_i=(\sqrt{2}-1)\,\goth i,
\mytag{3.15}\\
&\hskip -2em
C_i=-(\sqrt{2}+1)\,\goth i,
&&C_i=-(\sqrt{2}-1)\,\goth i.
\mytag{3.16}
\endxalignat
$$
Here $\goth i=\sqrt{-1}$. The roots \mythetag{3.16} are excluded by the
condition \mythetag{2.1}. The remain is two root \mythetag{3.15} of 
multiplicity $1$. They yield the following asymptotic expansions:
$$
\hskip -2em
\aligned
&t_i(p,q)=(\sqrt{2}+1)\,\goth i\,q^{\kern 0.7pt 2}\biggl(1+\sum^\infty_{s=1}\beta_{is}
\,q^{-s}\biggr),\\
&t_i(p,q)=(\sqrt{2}-1)\,\goth i\,q^{\kern 0.7pt 2}\biggl(1+\sum^\infty_{s=1}\beta_{is}
\,q^{-s}\biggr). 
\endaligned
\mytag{3.17}
$$
The results \mythetag{3.5}, \mythetag{3.11}, \mythetag{3.17} are summed up
in the following theorem. 
\mytheorem{3.1} For sufficiently large positive values of the parameter $q$, 
i\.\,e\. for $q>q_{\text{\,min}}$, the tenth-degree equation \mythetag{1.15} 
has five roots of multiplicity $1$ satisfying the condition \mythetag{2.1}. 
Three of them $t_1$, $t_2$, and $t_3$ are real roots. Their asymptotics as 
$q\to +\infty$ are given by the formulas 
$$
\xalignat 3
&\hskip -2em
t_1\sim p^{\kern 1pt 2},
&&t_2\sim p^{\kern 1pt 2},
&&t_3\sim p\,q.
\mytag{3.18}
\endxalignat
$$
The rest two roots $t_4$ and $t_5$ of the equation \mythetag{1.15} are complex. 
Their asymptotics as $q\to +\infty$ are given by the formulas 
$$
\xalignat 2
&\hskip -2em
t_4\sim (\sqrt{2}+1)\,\goth i\,q^{\kern 0.7pt 2},
&&t_5\sim (\sqrt{2}-1)\,\goth i\,q^{\kern 0.7pt 2}.
\mytag{3.19}
\endxalignat
$$
\endproclaim
The complex roots \mythetag{3.19} do not provide perfect cuboids. However, below
they are important for determining the exact number of real roots.\par
\head
4. Asymptotic estimates for real roots.
\endhead
     According to the formula \mythetag{3.18} the roots $t_1$ and $t_2$ are not
growing as $q\to +\infty$. For this reason we do not need to calculate 
$\beta_{is}$ in \mythetag{3.5} for them. But we need to find estimates 
for remainder terms $R_1$ and $R_2$ in the formulas
$$
\xalignat 2
&\hskip -2em
t_1=p^{\kern 1pt 2}+R_1(p,q), 
&&t_2=p^{\kern 1pt 2}+R_2(p,q) 
\mytag{4.1}
\endxalignat
$$
as $q\to +\infty$. Our goal is to obtain estimates of the form  
$$
\hskip -2em
|R_i(p,q)|<\frac{C(p)}{q}\text{, \ where \ }i=1,\,2. 
\mytag{4.2}
$$
In order to get such estimates we substitute 
$$
\hskip -2em
t=p^{\kern 1pt 2}+\frac{c}{q}
\mytag{4.3}
$$
into the equation \mythetag{1.15}. Then we perform another substitution into the 
equation obtained as a result of substituting \mythetag{4.3} into 
\mythetag{1.15}:
$$
\hskip -2em
q=\frac{1}{z}.
\mytag{4.4}
$$
Upon two substitutions \mythetag{4.3} and \mythetag{4.4} and upon removing 
denominators the equation \mythetag{1.15} is written as a polynomial equation 
in the new variables $c$ and $z$. It is a peculiarity of this equation that
it can be written as 
$$
\hskip -2em
16\,p^{\kern 1pt 12}+f(c,p,z)=4\,p^{\kern 1pt 6}\,c^2.
\mytag{4.5}
$$
Here $f(c,p,z)$ is a polynomial given by an explicit formula. The formula for
$f(c,p,z)$ is rather huge. Therefore it is placed to the ancillary file 
\darkred{{\tt strategy\kern -0.5pt\_\kern 1.5pt formulas.txt}} in a 
machine-readable form.\par
     Let $q\geqslant 59\,p$ and let the parameter $c$ run over the interval
from  $-5\,p^{\kern 1pt 3}$ to $0$:
$$
\hskip -2em
-5\,p^{\kern 1pt 3}<c<0.
\mytag{4.6}
$$
From $q\geqslant 59\,p$ and from \mythetag{4.4} we derive the estimate 
$|z|\leqslant 1/59\,p^{\kern 1pt -1}$. Using this estimate and using the 
inequalities \mythetag{4.6}, by means of direct calculations one can derive 
the following estimate for the modulus of the function $f(c,p,z)$:
$$
\hskip -2em
|f(c,p,z)|<15\,p^{\kern 1pt 12}.
\mytag{4.7}
$$
For fixed $p$ and $z$ the estimate \mythetag{4.7} means that the left hand side of 
the equation \mythetag{4.5} is a continuous function of $c$ taking the values 
within the range from $p^{\kern 1pt 12}$ to $31\,p^{\kern 1pt 12}$ while $c$ is 
in the interval \mythetag{4.6}. The right hand side of \mythetag{4.5} is also a 
continuous function of $c$. It decreases from $100\,p^{12}$ to $0$ in the interval \mythetag{4.6}. Therefore somewhere in the interval \mythetag{4.6} there is at 
least one root of the equation \mythetag{4.5}.\par 
     The parameter $c$ is related to the initial variable $t$ by means of the
formula \mythetag{4.3}. The inequalities \mythetag{4.5} for $c$ imply the
following inequalities for $t$:
$$
\hskip -2em
p^{\kern 1pt 2}-\frac{5\,p^{\kern 1pt 3}}{q}<t<p^{\kern 1pt 2}.
\mytag{4.8}
$$
The inequalities \mythetag{4.8} and the above considerations prove the following
theorem.
\mytheorem{4.1} For each $q\geqslant 59\,p$ there is at least one real root of 
the equation \mythetag{1.15} satisfying the inequalities \mythetag{4.8}.
\endproclaim
     The above considerations can be repeated for the case where the parameter 
$c$ runs over the interval from $0$ to $5\,p^{\kern 1pt 3}$. In this case due to 
\mythetag{4.3} from 
$$
0<c<5\,p^{\kern 1pt 3}
$$
we derive the inequalities 
$$
\hskip -2em
p^{\kern 1pt 2}<t<p^{\kern 1pt 2}+\frac{5\,p^{\kern 1pt 3}}{q}
\mytag{4.9}
$$
for the variable $t$ and hence we obtain the following theorem.\par 
\mytheorem{4.2} For each $q\geqslant 59\,p$ there is at least one real root of 
the equation \mythetag{1.15} satisfying the inequalities \mythetag{4.9}.
\endproclaim
     Now let's proceed to the growing root $t_3$ of the equation \mythetag{1.15} 
(see Theorem~\mythetheorem{3.1}). Upon refining the asymptotic formula  
\mythetag{3.18} for $t_3$ looks like 
$$
\hskip -2em
t_3=p\,q-\frac{16\,p^{\kern 1pt 3}}{q}+R_3(p,q).
\mytag{4.10}
$$
The formula \mythetag{4.10} is in agreement with the expansion 
\mythetag{3.11}. It means that 
$$
\xalignat 2
&\beta_{31}=0,
&&\beta_{32}=-16\,p^{\kern 1pt 2}.
\endxalignat
$$
Like in \mythetag{4.2}, our goal here is to obtain estimates of the form 
$$
\hskip -2em
|R_3(q)|<\frac{C(p)}{q^{\kern 0.7pt 2}}.
\mytag{4.11}
$$
In order to get such estimates we substitute 
$$
\hskip -2em
t=p\,q-\frac{16\,p^{\kern 1pt 3}}{q}+\frac{c}{q^{\kern 0.7pt 2}}
\mytag{4.12}
$$
into the equation \mythetag{1.15}. Immediately after that we perform the
substitution \mythetag{4.4} into the equation obtained by
substituting \mythetag{4.12} into \mythetag{1.15}. As a result of two 
substitutions \mythetag{4.12} and \mythetag{4.4} upon eliminating 
denominators the equation \mythetag{1.15} is written as a polynomial 
equation in the new variables $c$ and $z$. It looks like
$$
\hskip -2em
\varphi(c,p,z)=-2\,p^{\kern 1pt 5}\,c.
\mytag{4.13}
$$
Here $\varphi(c,p,z)$ is a polynomial of three variables. The explicit formula 
for $\varphi(c,p,z)$ is rather huge. Therefore it is placed to the ancillary 
file \darkred{{\tt strategy\kern -0.5pt\_\kern 1.5pt formulas.txt}} in a 
machine-readable form.\par
     Let $q\geqslant 59\,p$ and let the parameter $c$ run over the interval
from  $-5\,p^{\kern 1pt 4}$ to $5\,p^{\kern 1pt 4}$:
$$
\hskip -2em
-5\,p^{\kern 1pt 4}<c<5\,p^{\kern 1pt 4}.
\mytag{4.14}
$$
From $q\geqslant 59\,p$ and from \mythetag{4.4} we derive the estimate 
$|z|\leqslant 1/59\,p^{\kern 1pt -1}$. Using this estimate and using the 
inequalities \mythetag{4.14}, by means of direct calculations one can derive 
the following estimate for the modulus of the function $\varphi(c,p,z)$:
$$
\hskip -2em
|\varphi(c,z)|<7\,p^{\kern 1pt 9}.
\mytag{4.15}
$$
For fixed $p$ and $z$ the estimate \mythetag{4.15} means that the left hand 
side of the equation \mythetag{4.13} is a continuous function of $c$ taking 
the values within the range from $-7\,p^{\kern 1pt 9}$ to $7\,p^{\kern 1pt 9}$ 
while $c$ is in the interval \mythetag{4.14}. The right hand side of the equation 
\mythetag{4.13} is also a continuous function of $c$. It decreases from 
$10\,p^{\kern 1pt 9}$ to $-10\,p^{\kern 1pt 9}$ in the interval \mythetag{4.14}. 
Therefore somewhere in the interval \mythetag{4.14} there is at least one root 
of the polynomial equation \mythetag{4.13}.\par 
     The parameter $c$ is related to the initial variable $t$ by means of the
formula \mythetag{4.12}. The inequalities \mythetag{4.14} for $c$ imply the
following inequalities for $t$:
$$
\hskip -2em
p\,q-\frac{16\,p^{\kern 1pt 3}}{q}-\frac{5\,p^{\kern 1pt 4}}{q^{\kern 0.7pt 2}}
<t<p\,q-\frac{16\,p^{\kern 1pt 3}}{q}+\frac{5\,p^{\kern 1pt 4}}{q^{\kern 0.7pt 2}}.
\mytag{4.16}
$$
The inequalities \mythetag{4.16} and the above considerations prove the following
theorem.
\mytheorem{4.3} For each $q\geqslant 59\,p$ there is at least one real root of 
the equation \mythetag{1.15} satisfying the inequalities \mythetag{4.16}.
\endproclaim
    Theorems~\mythetheorem{4.1}, \mythetheorem{4.2}, and  \mythetheorem{4.3}
solve the problem of obtaining estimates of the form \mythetag{4.2} and
\mythetag{4.11} for the remainder terms in the refined asymptotic expansions
\mythetag{4.1} and \mythetag{4.10} for $q\geqslant 59\,p$.\par
\head
5. Asymptotic estimates for complex roots.
\endhead
     Let's proceed to complex roots of the equation \mythetag{1.15}. Upon 
refining the asymptotic formula \mythetag{3.19} for the complex root $t_4$ 
is written as 
$$
\hskip -2em
t_4=(\sqrt{2}+1)\,\goth i\,q^{\kern 0.7pt 2}+(\sqrt{2}-2)\,\goth i\,p^{\kern 1pt 2}+R_4(p,q)
\text{,\ \ where \ }\goth i=\sqrt{-1}.
\mytag{5.1}
$$
The formula \mythetag{5.1} is in agreement with the first expansion 
\mythetag{3.17}. It means that $\beta_{41}=0$ and 
$\beta_{42}=(4-3\,\sqrt{2})\,p^{\kern 1pt 2}$. \pagebreak Like in the formula 
\mythetag{4.2} and in the formula \mythetag{4.11}, our goal here is to 
obtain estimates of the form 
$$
\hskip -2em
|R_4(p,q)|<\frac{C(p)}{q}.
\mytag{5.2}
$$
In order to get such estimates we substitute 
$$
\hskip -2em
t=(\sqrt{2}+1)\,\goth i\,q^{\kern 0.7pt 2}+(\sqrt{2}-2)\,\goth i\,p^{\kern 1pt 2}
+\frac{c\,\goth i}{q}
\mytag{5.3}
$$
into the equation \mythetag{1.15}. Immediately after that we perform the
substitution \mythetag{4.4} into the equation obtained by substituting 
\mythetag{5.3} into \mythetag{1.15}. As a result of two 
substitutions \mythetag{5.3} and \mythetag{4.4} upon eliminating 
denominators the equation \mythetag{1.15} is written as a polynomial 
equation in the new variables $c$ and $z$. It looks like 
$$
\hskip -2em
\psi(z,p,c)=16\,c.
\mytag{5.4}
$$
Here $\psi(c,p,z)$ is a polynomial of three variables with purely real 
coefficients. The explicit formula for $\psi(c,p,z)$ is rather huge. Therefore 
it is placed to the ancillary file \darkred{{\tt strategy\kern -0.5pt\_\kern 1.5pt formulas.txt}} in a machine-readable form.\par
     Let $q\geqslant 59\,p$ and let the parameter $c$ run over the interval
from  $-5\,p^{\kern 1pt 3}$ to $5\,p^{\kern 1pt 3}$:
$$
\hskip -2em
-5\,p^{\kern 1pt 3}<c<5\,p^{\kern 1pt 3}.
\mytag{5.5}
$$
From $q\geqslant 59\,p$ and from \mythetag{4.4} we derive the estimate 
$|z|\leqslant 1/59\,p^{\kern 1pt -1}$. Using this estimate and using the 
inequalities \mythetag{5.5}, by means of direct calculations one can derive 
the following estimate for the modulus of the function $\psi(c,p,z)$:
$$
\hskip -2em
|\psi(c,p,z)|<15\,p^{\kern 1pt 3}.
\mytag{5.6}
$$
For fixed $p$ and $z$ the estimate \mythetag{5.6} means that the left hand 
side of the equation \mythetag{5.4} is a continuous function of $c$ taking 
its values within the range from $-15\,p^{\kern 1pt 3}$ to $15\,p^{\kern 1pt 3}$
while $c$ runs over the interval \mythetag{5.5}. The right hand side of the equation 
\mythetag{5.4} is also a continuous function of $c$. It increases from 
$-80\,p^{\kern 1pt 3}$ to $80\,p^{\kern 1pt 3}$ in the interval \mythetag{5.5}. 
Therefore somewhere in the interval \mythetag{5.5} there is at least one root 
of the polynomial equation \mythetag{5.4}.\par 
     The parameter $c$ is related to the initial variable $t$ by means of the
formula \mythetag{5.3}. Therefore the inequalities \mythetag{5.5} for $c$ 
imply the following inequalities for $t$:
$$
(\sqrt{2}+1)\,q^{\kern 0.7pt 2}+(\sqrt{2}-2)\,p^{\kern 1pt 2}-\frac{5\,p^{\kern 1pt 3}}{q}
<\Img\,t<(\sqrt{2}+1)\,q^{\kern 0.7pt 2}+(\sqrt{2}-2)\,p^{\kern 1pt 2}
+\frac{5\,p^{\kern 1pt 3}}{q}.
\quad
\mytag{5.7}
$$
The inequalities \mythetag{5.7} and the above considerations prove the following
theorem.
\mytheorem{5.1} For each $q\geqslant 59\,p$ there is at least one purely imaginary
root of the equation \mythetag{1.15} satisfying the inequalities \mythetag{5.7}.
\endproclaim
     The complex root $t_5$ is similar to the root $t_4$. Upon refining the 
asymptotic formula \mythetag{3.19} for the complex root $t_5$ is written as 
$$
\hskip -2em
t_4=(\sqrt{2}-1)\,\goth i\,q^{\kern 0.7pt 2}+(\sqrt{2}+2)\,\goth i\,p^{\kern 1pt 2}+R_5(p,q)
\text{,\ \ where \ }\goth i=\sqrt{-1}.
\mytag{5.8}
$$
The formula \mythetag{5.8} is in agreement with the second expansion 
\mythetag{3.17}. It means that $\beta_{51}=0$ and 
$\beta_{52}=(4+3\,\sqrt{2})\,p^{\kern 1pt 2}$. Like in the formulas 
\mythetag{4.2}, \mythetag{4.11}, and \mythetag{4.2}, our goal here is to 
obtain estimates of the form 
$$
\hskip -2em
|R_5(p,q)|<\frac{C(p)}{q}.
\mytag{5.9}
$$
In order to get such estimates we substitute 
$$
\hskip -2em
t=(\sqrt{2}-1)\,\goth i\,q^{\kern 0.7pt 2}+(\sqrt{2}+2)\,\goth i\,p^{\kern 1pt 2}
+\frac{c\,\goth i}{q}
\mytag{5.10}
$$
into the equation \mythetag{1.15}. Immediately after that we perform the
substitution \mythetag{4.4} into the equation obtained by substituting 
\mythetag{5.10} into \mythetag{1.15}. As a result of two 
substitutions \mythetag{5.10} and \mythetag{4.4} upon eliminating 
denominators the equation \mythetag{1.15} is written as a polynomial 
equation in the new variables $c$ and $z$. It looks like 
$$
\hskip -2em
\eta(z,c)=16\,c.
\mytag{5.11}
$$
Here $\eta(c,p,z)$ is a polynomial of three variables. The explicit formula 
for $\eta(c,p,z)$ is rather huge. Therefore it is placed to the ancillary 
file \darkred{{\tt strategy\kern -0.5pt\_\kern 1.5pt formulas.txt}} in a 
machine-readable form.\par
     Let $q\geqslant 59\,p$ and let the parameter $c$ run over the interval
from  $-5\,p^{\kern 1pt 3}$ to $5\,p^{\kern 1pt 3}$ (see \mythetag{5.5}). 
From $q\geqslant 59\,p$ and from \mythetag{4.4} we derive the estimate 
$|z|\leqslant 1/59\,p^{\kern 1pt -1}$. Using this estimate and using the 
inequalities \mythetag{5.5}, by means of direct calculations one can derive 
the following estimate for the modulus of the function $\eta(c,p,z)$:
$$
\hskip -2em
|\eta(c,p,z)|<15\,p^{\kern 1pt 3}.
\mytag{5.12}
$$
For fixed $p$ and $z$ the estimate \mythetag{5.12} means that the left hand 
side of the equation \mythetag{5.11} is a continuous function of $c$ taking 
its values within the range from $-15\,p^{\kern 1pt 3}$ to $15\,p^{\kern 1pt 3}$
while $c$ runs over the interval \mythetag{5.5}. The right hand side of the 
equation \mythetag{5.12} is also a continuous function of $c$. It increases 
from $-80\,p^{\kern 1pt 3}$ to $80\,p^{\kern 1pt 3}$ in the interval 
\mythetag{5.5}. Therefore somewhere in the interval \mythetag{5.5} there is 
at least one root of the polynomial equation \mythetag{5.12}.\par 
     The parameter $c$ is related to the initial variable $t$ by means of the
formula \mythetag{5.10}. Therefore the inequalities \mythetag{5.5} for $c$ 
imply the following inequalities for $t$:
$$
(\sqrt{2}-1)\,q^{\kern 0.7pt 2}+(\sqrt{2}+2)\,p^{\kern 1pt 2}-\frac{5\,p^{\kern 1pt 3}}{q}
<\Img\,t<(\sqrt{2}-1)\,q^{\kern 0.7pt 2}+(\sqrt{2}+2)\,p^{\kern 1pt 2}
+\frac{5\,p^{\kern 1pt 3}}{q}.
\quad
\mytag{5.13}
$$
The inequalities \mythetag{5.13} and the above considerations prove the following
theorem.
\mytheorem{5.2} For each $q\geqslant 59\,p$ there is at least one purely imaginary
root of the equation \mythetag{1.15} satisfying the inequalities \mythetag{5.13}.
\endproclaim
    Theorems~\mythetheorem{5.1} and \mythetheorem{5.2} solve the problem of 
obtaining estimates of the form \mythetag{5.2} and
\mythetag{5.9} for the remainder terms in the refined asymptotic expansions
\mythetag{5.1} and \mythetag{5.8} for $q\geqslant 59\,p$. Along with 
Theorems~\mythetheorem{4.1}, \mythetheorem{4.2}, and \mythetheorem{4.3}, they
separate the roots $t_1$, $t_2$, $t_3$, $t_4$, $t_5$ of the equation 
\mythetag{1.15} from each other for sufficiently large $q$ and provide rather 
precise intervals for their location.\par
\head
6. Non-intersection of asymptotic intervals.
\endhead
     Theorems~\mythetheorem{4.1}, \mythetheorem{4.2}, \mythetheorem{4.3},
\mythetheorem{5.1}, \mythetheorem{5.2} define five asymptotic intervals 
\mythetag{4.8}, \mythetag{4.9}, \mythetag{4.16}, \mythetag{5.7}, and
\mythetag{5.13} for $q\geqslant 59\,p$. It is easy to see that the intervals
\mythetag{4.8} and \mythetag{4.9} do not intersect. For the other pairs of 
intervals among \mythetag{4.8}, \mythetag{4.9}, \mythetag{4.16}, \mythetag{5.7}, 
\mythetag{5.13} this is not so obvious. Therefore we need some elementary 
lemmas.
\mylemma{6.1} For $q\geqslant 59\,p$ the asymptotic intervals 
\mythetag{4.8}, \mythetag{4.9}, \mythetag{4.16}, \mythetag{5.7}, and
\mythetag{5.13} do not comprise the origin. 
\endproclaim
\demo{Proof} Indeed, from $q\geqslant 59\,p$ for the left endpoint of the 
interval \mythetag{4.8} we derive
$$
\hskip -2em
p^{\kern 1pt 2}-\frac{5\,p^{\kern 1pt 3}}{q}\geqslant p^{\kern 1pt 2}
-\frac{5\,p^{\kern 1pt 2}}{59}=\frac{54\,p^{\kern 1pt 2}}{59}>0.
\mytag{6.1}
$$ 
The left endpoint of the interval \mythetag{4.9} is obviously positive:
$p^{\kern 1pt 2}>0$. In the case of the interval \mythetag{4.16} from 
$q\geqslant 59\,p$ we derive 
$$
\hskip -2em
p\,q-\frac{16\,p^{\kern 1pt 3}}{q}-\frac{5\,p^{\kern 1pt 4}}{q^{\kern 0.7pt 2}}\geqslant
\frac{204430\,p^{\kern 1pt 2}}{3481}>58\,p^{\kern 1pt 2}>0.
\mytag{6.2}
$$
In the case of the imaginary intervals \mythetag{5.7} and \mythetag{5.13}
from $q\geqslant 59\,p$ we derive
$$
\hskip -2em
\aligned
&(\sqrt{2}+1)\,q^{\kern 0.7pt 2}+(\sqrt{2}-2)\,p^{\kern 1pt 2}-\frac{5\,p^{\kern 1pt 3}}{q}
>8403\,p^{\kern 1pt 2}>0,\\
&(\sqrt{2}-1)\,q^{\kern 0.7pt 2}+(\sqrt{2}+2)\,p^{\kern 1pt 2}-\frac{5\,p^{\kern 1pt 3}}{q}
>1445\,p^{\kern 1pt 2}>0.
\endaligned
\mytag{6.3}
$$
The above inequalities \mythetag{6.1}, \mythetag{6.2}, and
\mythetag{6.3} prove Lemma~\mythelemma{6.1}.
\qed\enddemo
     Lemma~\mythelemma{6.1} means that for $q\geqslant 59\,p$ the real 
intervals \mythetag{4.8}, \mythetag{4.9}, and \mythetag{4.16} do not 
intersect with the purely imaginary intervals \mythetag{5.7} and 
\mythetag{5.13}. Moreover, the inequalities \mythetag{6.1}, \mythetag{6.2}, 
and \mythetag{6.3} show that all of these intervals are located within 
positive half-lines of the real and imaginary axes. Therefore any roots of 
the equation  \mythetag{1.15} enclosed within these intervals satisfy
the condition \mythetag{2.1}. 
\mylemma{6.2} For $q\geqslant 59\,p$ the real asymptotic intervals 
\mythetag{4.8}, \mythetag{4.9}, and \mythetag{4.16} do not intersect 
with each other. 
\endproclaim
\demo{Proof} The open intervals \mythetag{4.8} and \mythetag{4.9} are 
adjacent. They have one common endpoint $t=p^{\kern 1pt 2}$, but this 
endpoint does not belong to them. Therefore the intervals \mythetag{4.8} 
and \mythetag{4.9} do not intersect with each other.\par
     In order to prove Lemma~\mythelemma{6.2} it is sufficient to
compare the right endpoint of the interval \mythetag{4.9} with the
left endpoint of the interval \mythetag{4.16}. From $q\geqslant 59\,p$
we derive
$$
\hskip -2em
p^{\kern 1pt 2}+\frac{5\,p^{\kern 1pt 3}}{q}
\leqslant\frac{64\,p^{\kern 1pt 2}}{59}<2\,p^{\kern 1pt 2}.
\mytag{6.4}
$$
Comparing \mythetag{6.4} with \mythetag{6.2}, we conclude that
$$
\pagebreak
\hskip -2em
p^{\kern 1pt 2}+\frac{5\,p^{\kern 1pt 3}}{q}<p\,q
-\frac{16\,p^{\kern 1pt 3}}{q}-\frac{5\,p^{\kern 1pt 4}}{q^{\kern 0.7pt 2}}
\mytag{6.5}
$$
for $q\geqslant 59\,p$. The inequality \mythetag{6.5} completes
the proof of Lemma~~\mythelemma{6.2}.
\qed\enddemo
\mylemma{6.3} For $q\geqslant 59\,p$ the imaginary asymptotic intervals 
\mythetag{5.7} and \mythetag{5.13} do not intersect with each other. 
\endproclaim
\demo{Proof} In order to prove Lemma~\mythelemma{6.3} it is sufficient 
to compare the bottom endpoint of the interval \mythetag{5.7} with the
top endpoint of the interval \mythetag{5.13}: 
$$
\hskip -2em
\Img\bigl(t_{\sssize\text{bottom}}^{\sssize\text{(5.7)}}-
t_{\sssize\text{top}}^{\sssize\text{(5.13)}}\bigr)
=2\,q^{\kern 0.7pt 2}-4\,p^{\kern 1pt 2}-\frac{10\,p^{\kern 1pt 3}}{q}.
\mytag{6.6}
$$
From $q\geqslant 59\,p$ and from \mythetag{6.6} we derive the inequalities
$$
\hskip -2em
\Img\bigl(t_{\sssize\text{bottom}}^{\sssize\text{(5.7)}}-
t_{\sssize\text{top}}^{\sssize\text{(5.13)}}\bigr)
\geqslant\frac{410512\,p^{\kern 1pt 2}}{59}>6957\,p^{\kern 1pt 2}>0.
\mytag{6.7}
$$
The inequalities \mythetag{6.7} complete the proof of Lemma~\mythelemma{6.3}.
\qed\enddemo
     Lemmas~\mythelemma{6.1}, \mythelemma{6.2}, and \mythelemma{6.3} are 
summed up in the following theorem.
\mytheorem{6.1} For $q\geqslant 59\,p$ five roots $t_1$, $t_2$, $t_3$, 
$t_4$, $t_5$ of the equation \mythetag{1.15} obeying the condition 
\mythetag{2.1} are simple. They are located within five disjoint intervals 
\mythetag{4.8}, \mythetag{4.9}, \mythetag{4.16}, \mythetag{5.7}, \mythetag{5.13}, 
one per each interval. 
\endproclaim
     Due to \mythetag{2.2} Theorem~\mythetheorem{6.1} locates all of the ten
roots of the equation \mythetag{1.15}. 
\head
7. Integer points of asymptotic intervals.
\endhead
     It is easy to see that the asymptotic intervals \mythetag{4.8}, \mythetag{4.9}, \mythetag{4.16} become more and more narrow if $p$ is fixed and $q\to+\infty$. 
Using this observation, one can easily prove the following two theorems. 
\mytheorem{7.1} If $q\geqslant 59\,p$ and $q>5\,p^{\kern 1pt 3}$, then the 
asymptotic intervals \mythetag{4.8} and \mythetag{4.9} have no integer points. 
\endproclaim
\mytheorem{7.2} If $q\geqslant 59\,p$ and $q^{\kern 0.7pt 2}>10\,p^{\kern 1pt 4}$, then the 
asymptotic interval \mythetag{4.16} has at most one integer point. 
\endproclaim
     The next theorem is more complicated. 
\mytheorem{7.3} If $q\geqslant 59\,p$ and $q\geqslant 16\,p^{\kern 1pt 3}
+5\,p/16$, then the asymptotic interval \mythetag{4.16} has no integer points. 
\endproclaim
\demo{Proof} Note that $(1+x)^2>1+2\,x$ for any positive $x$. This inequality
yields 
$$
\hskip -2em
1+x>\sqrt{1+2\,x}\text{\ \ for any \ }x>0.
\mytag{7.1}
$$
Let's write the inequality $q\geqslant 16\,p^{\kern 1pt 3}+5\,p/16$ in the 
following way: 
$$
\hskip -2em
q-8\,p^{\kern 1pt 3}\geqslant 8\,p^{\kern 1pt 3}+\frac{5\,p}{16}
=8\,p^{\kern 1pt 3}\,\Bigl(1+\frac{5}{128\,p^{\kern 1pt 2}}\Bigr).
\mytag{7.2}
$$
Setting $x=5/(128\,p^{\kern 1pt 2})$ in \mythetag{7.1} and applying it to 
\mythetag{7.2}, we get
$$
\hskip -2em
q-8\,p^{\kern 1pt 3}>8\,p^{\kern 1pt 3}\,\sqrt{1+\frac{5}{64\,p^{\kern 1pt 2}}}.
\mytag{7.3}
$$
Both sides of the inequality \mythetag{7.3} are positive. Squaring them, we obtain
$$
\hskip -2em
(q-8\,p^{\kern 1pt 3})^2>64\,p^{\kern 1pt 6}\,\Bigl(1
+\frac{5}{64\,p^{\kern 1pt 2}}\Bigr).
\mytag{7.4}
$$
Expanding both sides of the inequality \mythetag{7.4}, we bring it to
$$
\hskip -2em
q^{\kern 0.7pt 2}-16\,p^{\kern 1pt 3}\,q>5\,p^{\kern 1pt 4}.
\mytag{7.5}
$$
And finally, dividing both sides of the inequality \mythetag{7.5} by $q^{\kern 0.7pt 2}$, we 
write it as 
$$
\hskip -2em
\frac{16\,p^{\kern 1pt 3}}{q}+\frac{5\,p^{\kern 1pt 4}}{q^{\kern 0.7pt 2}}<1.
\mytag{7.6}
$$
\par
      Apart from $q\geqslant 16\,p^{\kern 1pt 3}+5\,p/16$, we have the inequality
$q\geqslant 59\,p$ that yields 
$$
\hskip -2em
0<\frac{5\,p^{\kern 1pt 4}}{q^{\kern 0.7pt 2}}\leqslant\frac{5\,p^{\kern 1pt 4}}{59\,p\,q}
=\frac{5\,p^{\kern 1pt 3}}{59\,q}<\frac{16\,p^{\kern 1pt 3}}{q}.
\mytag{7.7}
$$
Applying \mythetag{7.6} and \mythetag{7.7} to \mythetag{4.16}, 
we derive the following inequalities:
$$
\hskip -2em
p\,q-1<t<p\,q.
\mytag{7.8}
$$
Since $p\,q$ is integer, the inequalities \mythetag{7.8} have no integer
solutions for $t$. This means that the interval \mythetag{4.16} has no integer
points. Theorem~\mythetheorem{7.3} is proved. 
\qed\enddemo
\head
8. A strategy for numeric search.
\endhead
     The numeric search for perfect cuboids in the case of the second cuboid
conjecture (see Conjecture~\mytheconjecture{1.1}) is based on the equation 
\mythetag{1.15}. The equation \mythetag{1.15} is related to the equation 
\mythetag{1.8} through the substitutions \mythetag{1.10} and \mythetag{1.13}. 
Substituting either \mythetag{1.10} or \mythetag{1.13} into the inequalities 
$t>a$, $t>b$, and $t>u$ from Theorem~\mythetheorem{1.2}, we get the same result 
expressed by the inequalities 
$$
\xalignat 3
&\hskip -2em
t>p^{\kern 1pt 2},
&&t>p\,q,
&&t>q^{\kern 0.7pt 2}.
\mytag{8.1}
\endxalignat
$$
Similarly, substituting either \mythetag{1.10} or \mythetag{1.13} into the 
inequality $(a+t)\,(b+t)>2\,t^2$, we get the same result expressed by
the inequality
$$
	\hskip -2em
(p^{\kern 1pt 2}+t)\,(p\,q+t)>2\,t^2.
\mytag{8.2}
$$
Theorem~\mythetheorem{1.2} specified for the case of second cuboid conjecture
(see Conjecture~\mytheconjecture{1.1}) is formulated in the following way.
\mytheorem{8.1} A triple of integer numbers $p$, $q$, and $t$ satisfying the
equation \mythetag{1.15} and such that $p\neq q$ are coprime provides a perfect
cuboid of and only if the inequalities \mythetag{8.1} and \mythetag{8.2} are
fulfilled. 
\endproclaim
     Generally speaking, the numeric search based on Theorem~\mythetheorem{8.1} 
is a three-para\-metric search. The inequalities \mythetag{8.1} set lower bounds 
for $t$, but they do not restrict $t$ to a finite set of values. The inequality 
\mythetag{8.2} is different. Since $p$ and $q$ are positive, the inequality 
\mythetag{8.2} can be written as follows:
$$
\hskip -2em
t<\frac{p^{\kern 1pt 2}+p\,q}{2}+\frac{p\,\sqrt{p^{\kern 1pt 2}+6\,p\,q
+q^{\kern 0.7pt 2}}}{2}.
\mytag{8.3}
$$
Due to \mythetag{8.1} and \mythetag{8.3} for each fixed $p$ and fixed $q$ one 
should iterate only over a finite set of integer values of $t$, i\.\,e\. the search 
is two-parametric in effect.\par
     Theorems~\mythetheorem{7.1} and \mythetheorem{7.3} strengthen the restrictions.
They say that for each fixed $p$ one should iterate over a finite set of values 
of $q$ and $t$, i\.\,e\. the search for perfect cuboids becomes effectively 
one-parametric.\par
     Theorems~\mythetheorem{7.1} and \mythetheorem{7.3} can be further strengthened
with the use of the inequalities \mythetag{8.1}. Assume that the condition
$q\geqslant 59\,p$ is fulfilled and assume that $t$ belongs to the first
asymptotic interval \mythetag{4.8} (see Theorem~\mythetheorem{6.1}). Then from
\mythetag{8.1} and \mythetag{4.8} we derive two contradictory inequalities 
$t>p^{\kern 1pt 2}$ and $t<p^{\kern 1pt 2}$. Hence we have the following theorem.
\mytheorem{8.2} If $q\geqslant 59\,p$, then the asymptotic interval \mythetag{4.8} 
has no points satisfying the inequalities \mythetag{8.1}.
\endproclaim
     Similarly, if $q\geqslant 59\,p$ and if $t$ belongs to the second asymptotic 
interval \mythetag{4.9}, then from the inequalities \mythetag{4.9} and \mythetag{8.1} 
we derive
$$
\hskip -2em
\aligned
&t>q^{\kern 0.7pt 2}\geqslant (59\,p)^2=3481\,p^{\kern 1pt 2},\\
&t<p^{\kern 1pt 2}+\frac{5\,p^{\kern 1pt 3}}{q}.
\endaligned
\mytag{8.4}
$$
The inequalities \mythetag{8.4} mean that 
$$
\hskip -2em
3481\,p^{\kern 1pt 2}<p^{\kern 1pt 2}+\frac{5\,p^{\kern 1pt 3}}{q}.
\mytag{8.5}
$$
From \mythetag{8.5} one easily derives the inequality $q<5\,p/3480=p/696$ that
contradicts $q\geqslant 59\,p$. Hence we have the following theorem.
\mytheorem{8.3} If $q\geqslant 59\,p$, then the asymptotic interval \mythetag{4.9} 
has no points satisfying the inequalities \mythetag{8.1}.
\endproclaim
     Finally, assume that $q\geqslant 59\,p$ and let $t$ belong to the
third asymptotic interval \mythetag{4.16}. We know that $q\geqslant 59\,p$
implies \mythetag{7.7}. From \mythetag{7.7} and \mythetag{4.16} we derive
$$
\hskip -2em
t<p\,q-\frac{16\,p^{\kern 1pt 3}}{q}+\frac{5\,p^{\kern 1pt 4}}{q^{\kern 0.7pt 2}}<p\,q.
\mytag{8.6}
$$
The inequalities \mythetag{8.6} contradict the inequality $t>p\,q$ from 
\mythetag{8.1}. This contradiction proves the following theorem. 
\mytheorem{8.4} If $q\geqslant 59\,p$, then the asymptotic interval \mythetag{4.16} 
has no points satisfying the inequalities \mythetag{8.1}.
\endproclaim
     The fourth and fifth asymptotic intervals \mythetag{5.7} and \mythetag{5.13}
are purely imaginary. They do not provide perfect cuboids. \pagebreak Therefore 
Theorems~\mythetheorem{8.1}, \mythetheorem{8.2}, \mythetheorem{8.3}, and
\mythetheorem{8.4} are summarized in the following theorem.\par
\mytheorem{8.5} If $q\geqslant 59\,p$, then the Diophantine equation 
\mythetag{1.15} has no solutions providing perfect cuboids. 
\endproclaim
Theorem~\mythetheorem{8.5} is a background for an optimized strategy of 
numeric search for perfect cuboids. It says that only for $q<59\,p$ perfect 
cuboids are expected.\par 
     Applying the inequality $q<59\,p$\,\ to the formula \mythetag{8.3} we can 
simplify it as 
$$
\hskip -2em
t<61\,p^{\kern 1pt 2}.
\mytag{8.7}
$$
The upper bound \mythetag{8.7} is more simple than \mythetag{8.3}, though it
can be computationally less time-efficient for large values of $p$. Along with 
\mythetag{8.1} and the inequality 
$$
\hskip -2em
q<59\,p,
\mytag{8.8}
$$
it provides the following very simple computer code for our optimized strategy:
\medskip 
\def\myskip{\kern 1em}
\darkred{
\leftline{\tt\myskip for p from 1 by 1 to $+\infty$ do}
\leftline{\tt\myskip \ for q from 1 by 1 to 59*p-1 do}
\leftline{\tt\myskip \ \ if p<>q and gcd(p,q)=1}
\leftline{\tt\myskip \ \ \ then}
\leftline{\tt\myskip \ \ \ \ for t from max(p\shlyapka 2,p*q,q\shlyapka
2)+1 by 1 to 61*p\shlyapka2-1 do}
\leftline{\tt\myskip \ \ \ \ \ if Q\podcherkivanie pqt=0 and (p\shlyapka 
2+t)*(p*q+t)>2*t\shlyapka 2}
\leftline{\tt\myskip \ \ \ \ \ \ then}
\leftline{\tt\myskip \ \ \ \ \ \ \ Str:=sprintf("Cuboid is found: p=\%a,\ %
q=\%a, t=\%a.",p,q,t):}
\leftline{\tt\myskip \ \ \ \ \ \ \ writeline(default,Str):}
\leftline{\tt\myskip \ \ \ \ \ end if:}
\leftline{\tt\myskip \ \ \ \ end do:}
\leftline{\tt\myskip \ \ end if:}
\leftline{\tt\myskip \ end do:}
\leftline{\tt\myskip end do:}
}
\medskip\noindent
Here \darkred{{\tt gcd(p,q)}} stands for the greatest common divisor of $p$ 
and $q$, while \darkred{{\tt Q\podcherkivanie pqt}} is a computer version of 
the formula \mythetag{1.12}. In practice the infinity sign \darkred{$+\infty$}
should be replaced by some particular positive integer.\par
\head
9. Concluding remarks and acknowledgments. 
\endhead
     Though they look very simple, Theorem~\mythetheorem{8.5} and the 
inequality \mythetag{8.8} constitute the main result of the present work.
As for the above code, it should be further optimized e\.\,g\.
by some tricky algorithms for fast computing the values of 
\darkred{{\tt Q\podcherkivanie pqt}}. One of such further optimized 
versions of this code has been run on a desktop PC for $p$ from $1$ to 
$100$. It took about 7 hours for that. No perfect cuboids were found. 
\par
     The authors are grateful to Mr. Seth Kitchen from Missouri University 
of Science and Technology whose e-mail letter asking about a cuboid coding 
idea became an impetus for this work. 


\Refs
\ref\myrefno{1}\paper
\myhref{http://en.wikipedia.org/wiki/Euler\podcherkivanie 
brick}{Euler brick}\jour Wikipedia\publ 
Wikimedia Foundation Inc.\publaddr San Francisco, USA 
\endref
\ref\myrefno{2}\by Halcke~P.\book Deliciae mathematicae oder mathematisches 
Sinnen-Confect\publ N.~Sauer\publaddr Hamburg, Germany\yr 1719
\endref
\ref\myrefno{3}\by Saunderson~N.\book Elements of algebra, {\rm Vol. 2}\publ
Cambridge Univ\. Press\publaddr Cambridge\yr 1740 
\endref
\ref\myrefno{4}\by Euler~L.\book Vollst\"andige Anleitung zur Algebra, \rm
3 Theile\publ Kaiserliche Akademie der Wissenschaf\-ten\publaddr St\.~Petersburg
\yr 1770-1771
\endref
\ref\myrefno{5}\by Pocklington~H.~C.\paper Some Diophantine impossibilities
\jour Proc. Cambridge Phil\. Soc\. \vol 17\yr 1912\pages 108--121
\endref
\ref\myrefno{6}\by Dickson~L.~E\book History of the theory of numbers, 
{\rm Vol\. 2}: Diophantine analysis\publ Dover\publaddr New York\yr 2005
\endref
\ref\myrefno{7}\by Kraitchik~M.\paper On certain rational cuboids
\jour Scripta Math\.\vol 11\yr 1945\pages 317--326
\endref
\ref\myrefno{8}\by Kraitchik~M.\book Th\'eorie des Nombres,
{\rm Tome 3}, Analyse Diophantine et application aux cuboides 
rationelles \publ Gauthier-Villars\publaddr Paris\yr 1947
\endref
\ref\myrefno{9}\by Kraitchik~M.\paper Sur les cuboides rationelles
\jour Proc\. Int\. Congr\. Math\.\vol 2\yr 1954\publaddr Amsterdam
\pages 33--34
\endref
\ref\myrefno{10}\by Bromhead~T.~B.\paper On square sums of squares
\jour Math\. Gazette\vol 44\issue 349\yr 1960\pages 219--220
\endref
\ref\myrefno{11}\by Lal~M., Blundon~W.~J.\paper Solutions of the 
Diophantine equations $x^2+y^2=l^2$, $y^2+z^2=m^2$, $z^2+x^2
=n^2$\jour Math\. Comp\.\vol 20\yr 1966\pages 144--147
\endref
\ref\myrefno{12}\by Spohn~W.~G.\paper On the integral cuboid\jour Amer\. 
Math\. Monthly\vol 79\issue 1\pages 57-59\yr 1972 
\endref
\ref\myrefno{13}\by Spohn~W.~G.\paper On the derived cuboid\jour Canad\. 
Math\. Bull\.\vol 17\issue 4\pages 575-577\yr 1974
\endref
\ref\myrefno{14}\by Chein~E.~Z.\paper On the derived cuboid of an 
Eulerian triple\jour Canad\. Math\. Bull\.\vol 20\issue 4\yr 1977
\pages 509--510
\endref
\ref\myrefno{15}\by Leech~J.\paper The rational cuboid revisited
\jour Amer\. Math\. Monthly\vol 84\issue 7\pages 518--533\yr 1977
\moreref see also Erratum\jour Amer\. Math\. Monthly\vol 85\page 472
\yr 1978
\endref
\ref\myrefno{16}\by Leech~J.\paper Five tables relating to rational cuboids
\jour Math\. Comp\.\vol 32\yr 1978\pages 657--659
\endref
\ref\myrefno{17}\by Spohn~W.~G.\paper Table of integral cuboids and their 
generators\jour Math\. Comp\.\vol 33\yr 1979\pages 428--429
\endref
\ref\myrefno{18}\by Lagrange~J.\paper Sur le d\'eriv\'e du cuboide 
Eul\'erien\jour Canad\. Math\. Bull\.\vol 22\issue 2\yr 1979\pages 239--241
\endref
\ref\myrefno{19}\by Leech~J.\paper A remark on rational cuboids\jour Canad\. 
Math\. Bull\.\vol 24\issue 3\yr 1981\pages 377--378
\endref
\ref\myrefno{20}\by Korec~I.\paper Nonexistence of small perfect 
rational cuboid\jour Acta Math\. Univ\. Comen\.\vol 42/43\yr 1983
\pages 73--86
\endref
\ref\myrefno{21}\by Korec~I.\paper Nonexistence of small perfect 
rational cuboid II\jour Acta Math\. Univ\. Comen\.\vol 44/45\yr 1984
\pages 39--48
\endref
\ref\myrefno{22}\by Wells~D.~G.\book The Penguin dictionary of curious and 
interesting numbers\publ Penguin publishers\publaddr London\yr 1986
\endref
\ref\myrefno{23}\by Bremner~A., Guy~R.~K.\paper A dozen difficult Diophantine 
dilemmas\jour Amer\. Math\. Monthly\vol 95\issue 1\yr 1988\pages 31--36
\endref
\ref\myrefno{24}\by Bremner~A.\paper The rational cuboid and a quartic surface
\jour Rocky Mountain J\. Math\. \vol 18\issue 1\yr 1988\pages 105--121
\endref
\ref\myrefno{25}\by Colman~W.~J.~A.\paper On certain semiperfect cuboids\jour
Fibonacci Quart.\vol 26\issue 1\yr 1988\pages 54--57\moreref see also\nofrills 
\paper Some observations on the classical cuboid and its parametric solutions
\jour Fibonacci Quart\.\vol 26\issue 4\yr 1988\pages 338--343
\endref
\ref\myrefno{26}\by Korec~I.\paper Lower bounds for perfect rational cuboids 
\jour Math\. Slovaca\vol 42\issue 5\yr 1992\pages 565--582
\endref
\ref\myrefno{27}\by Guy~R.~K.\paper Is there a perfect cuboid? Four squares 
whose sums in pairs are square. Four squares whose differences are square 
\inbook Unsolved Problems in Number Theory, 2nd ed.\pages 173--181\yr 1994
\publ Springer-Verlag\publaddr New York 
\endref
\ref\myrefno{28}\by Rathbun~R.~L., Granlund~T.\paper The integer cuboid table 
with body, edge, and face type of solutions\jour Math\. Comp\.\vol 62\yr 1994
\pages 441--442
\endref
\ref\myrefno{29}\by Rathbun~R.~L., Granlund~T.\paper The classical rational 
cuboid table of Maurice Kraitchik\jour Math\. Comp\.\vol 62\yr 1994
\pages 442--443
\endref
\ref\myrefno{30}\by Peterson~B.~E., Jordan~J.~H.\paper Integer hexahedra equivalent 
to perfect boxes\jour Amer\. Math\. Monthly\vol 102\issue 1\yr 1995\pages 41--45
\endref
\ref\myrefno{31}\by Van Luijk~R.\book On perfect cuboids, \rm Doctoraalscriptie
\publ Mathematisch Instituut, Universiteit Utrecht\publaddr Utrecht\yr 2000
\endref
\ref\myrefno{32}\by Rathbun~R.~L.\paper The rational cuboid table of Maurice 
Kraitchik\jour e-print \myhref{http://arxiv.org/abs/math/0111229}{math.HO/0111229} 
in Electronic Archive \myEarXivlink
\endref
\ref\myrefno{33}\by Narumiya~N., Shiga~H.\paper On Certain Rational Cuboid Problems
\jour Nihonkai Math. Journal\vol 12\yr 2001\issue 1\pages 75--88
\endref
\ref\myrefno{34}\by Hartshorne~R., Van Luijk~R.\paper Non-Euclidean Pythagorean 
triples, a problem of Euler, and rational points on K3 surfaces\publ e-print 
\myhref{http://arxiv.org/abs/math/0606700}{math.NT/0606700} 
in Electronic Archive \myEarXivlink
\endref
\ref\myrefno{35}\by Waldschmidt~M.\paper Open diophantine problems\publ e-print 
\myhref{http://arxiv.org/abs/math/0312440}{math.NT/0312440} 
in Electronic Archive \myEarXivlink
\endref
\ref\myrefno{36}\by Ionascu~E.~J., Luca~F., Stanica~P.\paper Heron triangles 
with two fixed sides\publ e-print \myhref{http://arxiv.org/abs/math/0608185}
{math.NT/0608} \myhref{http://arxiv.org/abs/math/0608185}{185} in Electronic 
Archive \myEarXivlink
\endref
\ref\myrefno{37}\by Ortan~A., Quenneville-Belair~V.\paper Euler's brick
\jour Delta Epsilon, McGill Undergraduate Mathematics Journal\yr 2006\vol 1
\pages 30-33
\endref
\ref\myrefno{38}\by Knill~O.\paper Hunting for Perfect Euler Bricks\jour Harvard
College Math\. Review\yr 2008\vol 2\issue 2\page 102\moreref
see also \myhref{http://www.math.harvard.edu/\volna knill/various/eulercuboid/index.html}
{http:/\negskp/www.math.harvard.edu/\textvolna knill/various/eulercuboid/index.html}
\endref
\ref\myrefno{39}\by Sloan~N.~J.~A\paper Sequences 
\myhref{http://oeis.org/A031173}{A031173}, 
\myhref{http://oeis.org/A031174}{A031174}, and \myhref{http://oeis.org/A031175}
{A031175}\jour On-line encyclopedia of integer sequences\publ OEIS Foundation 
Inc.\publaddr Portland, USA
\endref
\ref\myrefno{40}\by Stoll~M., Testa~D.\paper The surface parametrizing cuboids
\jour e-print \myhref{http://arxiv.org/abs/1009.0388}{arXiv:1009.0388} 
in Electronic Archive \myEarXivlink
\endref
\ref\myrefno{41}\by Meskhishvili~M.\paper Perfect cuboid and congruent number 
equation solutions\jour e-print 
\myhref{http://arxiv.org/abs/1211.6548}{arXiv:1211} 
\myhref{http://arxiv.org/abs/1211.6548}{.6548} 
in Electronic Archive \myEarXivlink
\endref
\ref\myrefno{42}\by Meskhishvili~M.\paper Parametric solutions for a 
nearly-perfect cuboid \jour e-print 
\myhref{http://arxiv.org/abs/1211.6548}{arXiv:1502.02375} 
in Electronic Archive \myEarXivlink
\endref
\ref\myrefno{43}\by Sharipov~R.~A.\paper A note on a perfect Euler cuboid.
\jour e-print \myhref{http://arxiv.org/abs/1104.1716}{arXiv:1104.1716} 
in Electronic Archive \myEarXivlink
\endref
\ref\myrefno{44}\by Sharipov~R.~A.\paper Perfect cuboids and irreducible 
polynomials\jour Ufa Mathematical Journal\vol 4, \issue 1\yr 2012\pages 153--160
\moreref see also e-print \myhref{http://arxiv.org/abs/1108.5348}{arXiv:1108.5348} 
in Electronic Archive \myEarXivlink
\endref
\ref\myrefno{45}\by Sharipov~R.~A.\paper A note on the first cuboid conjecture
\jour e-print \myhref{http://arxiv.org/abs/1109.2534}{arXiv:1109.2534} 
in Electronic Archive \myEarXivlink
\endref
\ref\myrefno{46}\by Sharipov~R.~A.\paper A note on the second cuboid conjecture.
Part~\uppercase\expandafter{\romannumeral 1} 
\jour e-print \myhref{http://arxiv.org/abs/1201.1229}{arXiv:1201.1229} 
in Electronic Archive \myEarXivlink
\endref
\ref\myrefno{47}\by Sharipov~R.~A.\paper A note on the third cuboid conjecture.
Part~\uppercase\expandafter{\romannumeral 1} 
\jour e-print \myhref{http://arxiv.org/abs/1203.2567}{arXiv:1203.2567} 
in Electronic Archive \myEarXivlink
\endref
\ref\myrefno{48}\by Sharipov~R.~A.\paper Perfect cuboids and multisymmetric 
polynomials\jour e-print \myhref{http://arxiv.org/abs/1203.2567}
{arXiv:1205.3135} in Electronic Archive \myEarXivlink
\endref
\ref\myrefno{49}\by Sharipov~R.~A.\paper On an ideal of multisymmetric polynomials 
associated with perfect cuboids\jour e-print \myhref{http://arxiv.org/abs/1206.6769}
{arXiv:1206.6769} in Electronic Archive \myEarXivlink
\endref
\ref\myrefno{50}\by Sharipov~R.~A.\paper On the equivalence of cuboid equations and 
their factor equations\jour e-print \myhref{http://arxiv.org/abs/1207.2102}
{arXiv:1207.2102} in Electronic Archive \myEarXivlink
\endref
\ref\myrefno{51}\by Sharipov~R.~A.\paper A biquadratic Diophantine equation 
associated with perfect cuboids\jour e-print 
\myhref{http://arxiv.org/abs/1207.4081}{arXiv:1207.4081} in Electronic Archive 
\myEarXivlink
\endref
\ref\myrefno{52}\by Ramsden~J\.~R\.\paper A general rational solution of an equation 
associated with perfect cuboids\jour e-print \myhref{http://arxiv.org/abs/1207.5339}{arXiv:1207.5339} in Electronic Archive 
\myEarXivlink
\endref
\ref\myrefno{53}\by Ramsden~J\.~R\., Sharipov~R.~A.\paper Inverse problems 
associated with perfect cuboids\jour e-print
\myhref{http://arxiv.org/abs/1207.6764}{arXiv:1207.6764}
in Electronic Archive \myEarXivlink
\endref
\ref\myrefno{54}\by Sharipov~R.~A.\paper On a pair of cubic equations 
associated with perfect cuboids \jour e-print
\myhref{http://arxiv.org/abs/1208.0308}{arXiv:1208}
\myhref{http://arxiv.org/abs/1208.0308}{.0308} in Electronic Archive \myEarXivlink
\endref
\ref\myrefno{55}\by Sharipov~R.~A.\paper On two elliptic curves associated 
with perfect cuboids \jour e-print
\myhref{http://arxiv.org/abs/1208.1227}{arXiv:1208.1227} in Electronic 
Archive \myEarXivlink
\endref
\ref\myrefno{56}\by Ramsden~J\.~R\., Sharipov~R.~A. \paper On singularities 
of the inverse problems associated with perfect cuboids  \jour e-print
\myhref{http://arxiv.org/abs/1208.1859}{arXiv:1208.1859} in 
Archive \myEarXivlink
\endref
\ref\myrefno{57}\by Ramsden~J\.~R\., Sharipov~R.~A. \paper On two algebraic 
parametrizations for rational solutions of the cuboid equations \jour e-print
\myhref{http://arxiv.org/abs/1208.2587}{arXiv:1208.2587} in Electronic 
Archive \myEarXivlink
\endref
\ref\myrefno{58}\by Sharipov~R.~A.\paper A note on solutions of the cuboid 
factor equations  \jour e-print
\myhref{http://arxiv.org/abs/1209.0723}{arXiv:1209.0723}
 in Electronic Archive \myEarXivlink
\endref
\ref\myrefno{59}\by Sharipov~R.~A.\paper A note on rational and elliptic curves 
associated with the cuboid factor equations \jour e-print
\myhref{http://arxiv.org/abs/1209.5706}{arXiv:1209.5706} in Electronic 
Archive \myEarXivlink
\endref
\ref\myrefno{60}\by Ramsden~J\.~R\., Sharipov~R.~A. \paper Two and three descent 
for elliptic curves associated with perfect cuboids \jour e-print
\myhref{http://arxiv.org/abs/1303.0765}{arXiv:1303.0765} in 
Archive \myEarXivlink
\endref
\ref\myrefno{61}\by Kashchenko~I\.~S\.\book Asymptotic expansions for solution 
of equations\publ RIO YarGU\publaddr Yaroslavl\yr 2011\moreref
see \myhref{http://math.uniyar.ac.ru/math/system/files/Kaschenko\podcherkivanie I.S.\podcherkivanie Asimptoticheskoe\podcherkivanie Razlozhenie.pdf}{http:/\negskp/math.uniyar.ac.ru/math/system/files/Kaschenko\_I.S.\_Asimptoticheskoe\_Raz}
\myhref{http://math.uniyar.ac.ru/math/system/files/Kaschenko\podcherkivanie I.S.\podcherkivanie Asimptoticheskoe\podcherkivanie Razlozhenie.pdf}{lozhenie.pdf}
\endref
\endRefs
\enddocument
\end